\documentclass[11pt,a4paper]{article}
\usepackage{caption2}

\usepackage{CJK}
\usepackage{tabls}
\usepackage{bm}
\usepackage{multirow}
\usepackage{amssymb}
\usepackage{amsfonts}
\usepackage{mathrsfs}
\usepackage{empheq}
\usepackage{amsmath}
\usepackage{amsthm}
\usepackage{graphicx}
\usepackage{color}
\usepackage{indentfirst}
\usepackage{hyperref}
\usepackage{float}
\usepackage{cases}
\usepackage[titletoc]{appendix}
\usepackage[colorinlistoftodos,english]{todonotes}

\allowdisplaybreaks%

\numberwithin{equation}{section}
\theoremstyle{definition}

\newtheorem{thm}{Theorem}[section]

\newtheorem{lem}{Lemma}[section]
\newtheorem{prop}{Proposition}[section]
\newtheorem{rem}{Remark}[section]
\newtheorem{defn}{Definition}[section]

\newcommand{\R}{\mathbb{R}}

\newcommand{\eps}{\varepsilon}
\newcommand{\supp}{\operatorname{supp}}

\newcommand{\la}{\lambda}
\newcommand{\ph}{\phi}
\newcommand{\pS}{p_{\mathrm S}}

\begin{document}

\begin{CJK*}{GB}{gbsn}
\title{Blow-up of damped semilinear wave equation with short range potential and critical power nonlinearity}

\author{Masahiro Ikeda \thanks{Graduate School of Information Science and Technology, The University of Osaka, 1-5 Yamadaoka, 565-0871, Suita, Osaka, Japan/Center for Advanced Intelligence Project RIKEN, 1-4-1, Nihonbashi, Chuo-ku, Tokyo 103-0027, Japan (ikeda@ist.osaka-u.ac.jp/masahiro.ikeda@a.riken.jp)}
\and Ziheng Tu \thanks{School of Data Science, Zhejiang University of Finance and Economics, 310018, Hangzhou, P.R.China (tuziheng@zufe.edu.cn)}
\and Chengbo Wang \thanks{School of Mathematical Sciences, Zhejiang Normal University, Jinhua 321004, China (wangcbo@zjnu.edu.cn)} }

\maketitle

\begin{abstract}
The present paper is devoted to the study of a semilinear wave equation with spatially dependent scattering damping and a short-range potential, with particular focus on determining the critical exponent for the power-type nonlinearity. By employing the test function method, we establish wave-type blow-up results and derive lifespan estimates at the critical level. The main contribution of this work lies in providing a uniform asymptotic analysis for a family of test functions, which is achieved through a detailed investigation of a parametric second-order elliptic equation. As a further novelty, we relax the radial symmetry condition conventionally imposed on the coefficients $D$ and $V$ by resorting to the comparison principle in conjunction with a monotone exhaustion argument.
\end{abstract}



\section{Introduction and main results}
The Cauchy problem for semilinear wave equations has been one of the central topics in nonlinear partial differential equations over the past several decades. Among the fundamental questions is the determination of the critical exponent separating global existence of small-data solutions from finite-time blow-up. For the prototype equation
\[
u_{tt}-\Delta u=|u|^p,
\qquad
(t,x)\in (0,\infty)\times\mathbb R^n,
\]
the behavior of solutions is governed by the Strauss exponent \(p_S(n)\), defined as the positive root of
\[
(n-1)p^2-(n+1)p-2=0.
\]
The Strauss conjecture predicts that sufficiently small initial data generate global solutions when
\(
p>p_S(n),
\)
whereas finite-time blow-up occurs whenever
\(
1<p\le p_S(n).
\)
This problem has stimulated a vast literature connecting harmonic analysis, nonlinear functional analysis, and the qualitative theory of nonlinear hyperbolic equations.

The study of blow-up phenomena began with the pioneering work of John \cite{John1979}, who proved finite-time blow-up for nonlinear wave equations in three space dimensions. Glassey \cite{Glassey1981} obtained further blow-up results by developing multiplier-type arguments. Sideris \cite{Sideris1984} introduced a robust concavity method for proving nonexistence of global solutions in high dimensions, while Schaeffer \cite{Schaeffer1985} treated the critical case for the equation
\[
u_{tt}-\Delta u=|u|^p.
\]
On the global existence side, the works of Lindblad and Sogge \cite{LindbladSogge1996}, Georgiev, Lindblad and Sogge \cite{GeorgievLindbladSogge1997}, and Tataru \cite{Tataru2000} developed weighted Strichartz estimates, local energy estimates, and refined space-time estimates. On the blow-up side, Yordanov and Zhang \cite{YordanovZhang2006} proved finite-time blow-up in the critical case together with sharp lifespan estimates. These results completed the classical picture for the constant-coefficient equation in the Strauss range.

Motivated by applications to wave propagation in inhomogeneous media, scattering theory, and mathematical physics, considerable attention has recently been devoted to wave equations containing lower-order perturbations. Typical examples include equations with damping, potential
\[
u_{tt}-\Delta u+D(x)u_t+V(x)u=|u|^p,
\]
where both the damping coefficient \(D(x)\) and the potential \(V(x)\) depend on the spatial variable. Such equations arise naturally in heterogeneous media, linearization around stationary states, acoustic and electromagnetic scattering, and effective models on noncompact geometric backgrounds. From the viewpoint of scattering theory, the spatial decay of
\(
D(x),\ V(x),
\)
determines whether the perturbation is short-range or long-range, and hence plays a decisive role in the large-time behavior of solutions.

For wave equation with damping, a systematic theory has been developed in the last two decades. Nishihara \cite{Nishihara2003} established \(L^p\)-\(L^q\) estimates for classical damped wave equations and applied them to nonlinear problems. For time-dependent damping with typical coefficients as ${\mu}{(1+t)^{-\beta}}$, Wirth \cite{Wirth04,Wirth06,Wirth07} provided a systematic classification according to the decay rate. This classification, which distinguishes among the over-damping ($\beta<0$), effective($0\leq\beta<1$), scale-invariant($\beta=1$), and non-effective(scattering, $\beta>1$) regimes, has subsequently stimulated further in-depth research on the semilinear problem. For example, see \cite{IkedaWakasugi2020}, \cite{DAbbiccoLucenteReissig2013,IkedaSobajimaWakasugi2019}, \cite{IkedaTuWakasa22} for overdamping, effective and scattering cases, respectively, and in particular \cite{Wakasugi2014,DAbbiccoLucenteReissig2015,PalmieriReissig2019,PalmieriTu2019} for the scale-invariant case. These works have fully clarified the transition from the heat-like to wave-like critical behavior as the strength of damping varies.

The study of the spatial dependent damping can be trace back to Mochizuki\cite{Mochizuki1976} for $0\leq D(x)\leq \langle x\rangle^{-\beta}$ with $\beta>1$ which showed that the solution behaves asymptotically like that of the linear wave equation. In the case $D(x)\backsim\langle x\rangle^{-\beta},\ \beta\in[0,1)$, Wakasugi\cite{Wakasugi2014_spatial} and Sobajima\cite{Sobajima2023} showed that the solution behaves like that of weighted heat equation and exhibits diffusion phenomena. See also \cite{IkedhataTakeda2020, RaduTodorovaYordanov2009, RaduTodorovaYordanov2010, RaduTodorovaYordanov2016}. These studies of the linear problem indicate that the spatial damping should share a classification structure in $\beta$ similar to that of time-dependent damping, and was further confirmed from the aspect of critical exponent of corresponding semilinear problem. For $\beta\in[0,1)$, Ikehata, Todorova and Yordanov \cite{IkedhataTodorovaYordanov2009} showed that the critical power is the shifted Fujita exponent $$p_c(n) = p_G(n-\beta+1)= 1+\frac{2}{n-\beta}.$$
The scale-invariant case \(D(x)\sim \mu|x|^{-1}\) is especially delicate.
Li \cite{Li2013} obtained a Fujita type blow-up result when $p\leq p_G(n)$ with $\mu\geq n$. Ikeda and Sobajima \cite{IkedaSobajima2021} established a Strauss type blow-up result for $0\leq\mu<\frac{(n-1)^2}{n+1}$ and $n-1<p\leq p_S(n +\mu)$, showing that the critical exponent is shifted by the quantity of damping. See also Dai, Kubo and Sobajima \cite{DaiKuboSobajima2021}. The global existence for $p>p_S(3+\mu)$ has recently been established for $\mu\in (0, 1/7]$ and $n=3$, in the work of Adamson, Metcalfe and Wang \cite{EMW26-2}.
In contrast, when \(D(x)=O(|x|^{-\beta}),\ \beta>1,\) the damping is of scattering type, and one expects the critical exponent to remain \(p_S(n)\). This expectation has been confirmed for the blow-up part by Lai and Tu \cite{LaiTu2020} for $\beta>2,\ n\geq3,$ and further extended by Lai, Liu, Tu and Wang \cite{LaiLiuTuWang2023} for $\beta>1$ and $n\geq2$. On the other hand,
in a recent work  \cite{EMW26}, Adamson, Metcalfe and Wang prove the global existence for $p>p_S(3)$ when $\beta>1$ and $n=3$.

The presence of a spatially dependent potential introduces an additional difficulty. Unlike damping, the potential modifies the elliptic structure associated with the wave operator. Therefore the classical test function method, which relies on explicit positive solutions of the adjoint equation, cannot be applied directly. Lai, Liu, Tu and Wang \cite{LaiLiuTuWang2023} developed a powerful framework for semilinear wave equations with space-dependent damping and potential. The key idea was to construct test functions from family of positive solutions of the parameter-dependent elliptic equation
\begin{equation}\label{eq:test_para}
\Delta \phi_\lambda
=
\bigl(\lambda^2+\lambda D(x)+V(x)\bigr)\phi_\lambda.
\end{equation}
They obtained sharp lifespan estimates in the subcritical regime. However, the critical case in the presence of a nontrivial short-range potential requires finer uniform asymptotic information on \(\phi_\lambda\) as \(\la\to0\) and \(|x|\to\infty\), and this remained the main unresolved point.

The main purpose of the present paper is to establish the critical blow-up result for semilinear wave equations with space-dependent scattering damping and short-range potentials. Moreover, the technical radial symmetric assumption on either $D$ or $V$ is removed through introducing radial barrier solution, which is another novelty. We consider the initial value problem
\begin{align}\label{eq:PDE0}
u_{tt} - \Delta u + D(x) u_t + V(x) u &= |u|^p,
\qquad (t,x) \in (0,T) \times \mathbb{R}^n,\notag\\
u(0,x)& =\varepsilon \, f(x) \qquad x\ \in\ \mathbb{R}^n, \\
u_t(0,x)& =\varepsilon \, g(x) \qquad \, x\ \in\ \mathbb{R}^n, \notag
\end{align}
where \(D\) and \(V\)
are nonnegative measurable functions satisfying the decay bounds
\begin{equation}
\label{eq:DV-main-ass}
0\le D(x)\le\frac{\mu_1}{(1+|x|)^\beta},
\qquad
0\le V(x)\le\frac{\mu_2}{(1+|x|)^\gamma}
\end{equation}
with $\beta>1$ and $\gamma>2$. As a preliminary step, we first establish the requisite asymptotics for the radial solution $\tilde{\phi}_\lambda(r)$ of the elliptic problem \eqref{eq:test_para}, in which the coefficients $D$ and $V$ are replaced by their radially symmetric majorants
\[
D_*(r)=\mu_1(1+r)^{-\beta},
\qquad
V_*(r)=\mu_2(1+r)^{-\gamma}.\]
More precisely, it is demonstrated that the radial solution exhibits the same leading-order behavior as the free growing solution in the exterior region, namely
\[
\tilde{\phi}_\lambda(r)
\sim
(1+\lambda r)^{-\frac{n-1}{2}}e^{\lambda r}
\ ,\ r>\lambda^{-1}.
\]
In contrast, within the interior region the solution is characterized by the normalization of zero-energy solution
\(
{\phi_0(r)}/{\phi_0(1/\lambda)},
\)
where $\phi_0(r)$ denotes the corresponding zero-energy limiting solution by taking $\lambda=0$ in \eqref{eq:test_para}. This limiting profile satisfies
\[
\phi_0(r)\sim   \left\{\begin{array}{l}\ln r,\quad n=2,\\
1,\quad n\geq3,\end{array}\right.\quad r\rightarrow\infty.
\]
The logarithmic growth in two dimension underscores a fundamental qualitative difference from the higher-dimensional case and therefore a more delicate, dimension-specific treatment is introduced. These radial solutions $\{\tilde{\phi}_\lambda\}$ subsequently serve as barriers in the non-radial setting. The construction of a global non-radial solution to the original problem is carried out by solving Dirichlet problems on expanding balls; the maximum principle then furnishes pointwise bounds that are independent of the radius of exhaustion. Equipped with these elliptic estimates, we build the dual test function \(b_q(t,x)\), which satisfies the desired decay properties for the critical test function method. We then prove finite-time blow-up at the Strauss critical exponent \(p=p_S(n)\), together with the corresponding upper bound for the lifespan of small-data solutions. The argument shows that the critical blow-up mechanism is governed by the asymptotic structure of the associated elliptic operator rather than by the explicit form of the coefficients.

We now give the definition of energy weak solution and the life-span of blowup solution.
\begin{defn} Let $f\in H^1(\mathbb{R}^n)$ and $g\in L^2(\mathbb{R}^n)$.
We say that $u$ is an energy solution of \eqref{eq:PDE0} on $[0,T)$ if
$$u\in C([0,T),H^1(\mathbb{R}^n))\cap C^1([0,T),L^2(\mathbb{R}^n))\cap L^p_{loc}(\mathbb{R}^n\times[0,T))$$
satisfies
\begin{align}
&\int_{\mathbb{R}^n}u_t(t,x)\phi(t,x)\,dx-\int_{\mathbb{R}^n}u_t(0,x)\phi(0,x)\,dx - \int_0^t\int_{\mathbb{R}^n} u_t(s,x)\phi_t(s,x) \,dx\, ds\notag \\
& \ +\int_0^t\int_{\mathbb{R}^n}\nabla u(s,x)\cdot\nabla\phi(s,x)\, dx\, ds+\int_0^t\int_{\mathbb{R}^n}\left(D(x) u_t(s,x)+V(x)u(s,x)\right)\phi(s,x)\,dx\, ds  \notag \\
& =\int_0^t \int_{\mathbb{R}^n}|u(s,x)|^p\phi(s,x)\,dx \, ds  \label{def}
\end{align}
for any $\phi\in C_0^\infty([0,T)\times\mathbb{R}^n)$ and any $t\in [0,T)$.

The supremum of all such time of existence, \( T \), is called the lifespan of the solution of \eqref{eq:PDE0}, denoted by \( T_\varepsilon \).

\end{defn}

The main results are stated as follows.

\begin{thm}[Critical blow-up and lifespan]
\label{thm:main}
Let \(n\ge2\), let \(D,V\) satisfy
\eqref{eq:DV-main-ass}, and let
\(
p=\pS(n).
\)
Assume
\(
f\in H^1(\R^n),\ g\in L^2(\R^n),
\
f,g\ge0,\ (f,g)\not\equiv(0,0),
\)
and
\(
\supp f\cup\supp g\subset B_R
\)
for some \(R>0\).
Let \(u\) be an energy solution of \eqref{eq:PDE0} on
\([0,T_\eps)\) satisfying
\begin{equation}
\label{eq:finite-prop-main}
\supp u(t,\cdot)\subset B_{t+R}.
\end{equation}
Then
\(
T_\eps<\infty.
\)
Moreover, for all sufficiently small \(\eps>0\),
\begin{equation}
\label{eq:lifespan-main}
T_\eps
\le
\exp\!\left(
C\eps^{-p(p-1)}
\right),
\end{equation}
where \(C>0\) is independent of \(\eps\).
\end{thm}

The paper is organized as follows. Section~2 investigates the parameter-dependent elliptic equation with radial majorants $D_*$ and $V_*$, and establishes uniform asymptotic bounds for the radial solution in both exterior and interior regimes. Section~3 derives the non-radial solution asymptotics for the elliptic equation with non-radial $D$ and $V$, which in turn yields the desired estimates for the test functions \(b_q\). Section~4 presents the critical blow-up result and its sharp lifespan upper bound. The appendix gathers the standard facts concerning local existence and finite speed of propagation.

\section{Asymptotic of radial barrier solution $\tilde{\phi}_\lambda(r)$}
In this section, we consider equation \eqref{eq:test_para} with $D$ and $V$ replaced by their majorants $D_*$ and $V_*$.
We write equation \eqref{eq:test_para} in radially symmetric form
\begin{equation}\label{eq:test_radial}
\tilde{\phi}_{\lambda}'' + \frac{n-1}{r} \tilde{\phi}_{\lambda}' - \left( \lambda D_* + V_* \right) \tilde{\phi}_{\lambda} = \lambda^2 \tilde{\phi}_{\lambda},
\end{equation}
or in its self-adjoint form
\begin{equation}\label{eq:adj_phi}
\left(r^{n-1} \tilde{\phi}_{\lambda}'\right)' = r^{n-1} \left( \lambda^2 + \lambda D_* +V_* \right) \tilde{\phi}_{\lambda}.
\end{equation}
Following the standard strategy, we first analyze the equation outside the ball \(B_c\) of radius \(r_c =1/\lambda\). Next, we investigate the lower bound for \(\tilde{\phi}_\lambda\) inside \(B_c\), where a two-step approximation procedure is employed. For notational simplicity, we retain the symbols \(D\) and \(V\) for their majorants and, with a slight abuse of notation, denote the radial barrier solution by \(\phi_\lambda\) in the following subsection.

\subsection{Asymptotic of $\phi_\lambda(r)$ outside $B_c$}
We first present the asymptotic behavior outside the ball $B_c$.
\begin{lem}\label{exterior}
Let $n\geqslant 2$, there exists a growing \( C^2 \) solution $\phi_\lambda(r)$ of \eqref{eq:test_radial} satisfying
\begin{equation}\label{outside}
 \phi_\lambda(r)  \simeq (\lambda r)^{-\frac{n-1}{2}}e^{\lambda r}\quad \text{for $r\geq r_c$}
\end{equation}
with any given \( 0 < \lambda < 1 \).
\end{lem}
\begin{proof}
The result follows directly by applying Lemma 3.2 in \cite{LiuWang23}. Here we give a complete proof. Plugging $\phi_{\lambda}(r) = r^{-\frac{n-1}{2}} y(r)$ into \eqref{eq:test_radial}, we obtain
\begin{equation}
y'' - \lambda^2 y + G(r) y = 0 \label{eq:y}
\end{equation}
where
\[
G(r) = -\left( \frac{(n-1)(n-3)}{4r^2} +\la D(r)+ V(r) \right).
\]
On the boundary of the ball $B(r_c)$, we set $y_0=y(r_c)=\lambda^{-\frac{n-1}2}$ and $y_1=y'(r_c)$.
We first give the upper bound of quotient ratio $\frac{y_1}{y_0}$. From the relation
\[y' = r^{\frac{n-1}{2}} \phi_{\lambda}'(r) + \frac{n-1}{2} r^{\frac{n-3}{2}} \phi_{\lambda}(r)
=r^{\frac{n-1}{2}} \phi_{\lambda}(r) \left[ \frac{\phi_{\lambda}'(r)}{\phi_{\lambda}(r)} + \frac{n-1}{2} r^{-1} \right],\]
it follows that
\[\frac{y_1}{y_0}= \frac{\phi_{\lambda}'(\frac1\lambda)}{\phi_{\lambda}(\frac1\lambda)} + \frac{n-1}{2} \lambda .
\]
Integrating \eqref{eq:adj_phi} from \(0\) to \(R\) yields
\[
 R^{n-1} \phi_{\lambda}'(R) \leq \int_{0}^{R} r^{n-1} \left( \lambda^2 + \lambda \frac{\mu_1}{1+r^{\beta}} + \frac{\mu_2}{1+r^{\gamma}} \right) \phi_{\lambda}(r) \, dr
\]
which implies
\begin{align}
\left|\dfrac{\phi_{\lambda}'(R)}{\phi_{\lambda}(R)}\right|\leq& R^{1-n} \cdot \int_{0}^{R} r^{n-1} \left( \lambda^2 + \lambda \frac{\mu_1}{1+r^{\beta}} + \frac{\mu_2}{1+r^{\gamma}} \right)  \, dr\label{ineq:phi}\\
\leq& R^{1-n} \cdot \left[ \frac{\lambda^2 r^n}{n} \bigg|_0^R + \frac{\lambda \mu_1 r^{n-\beta}}{n-\beta} \bigg|_0^R + \frac{\mu_2 r^{n-\gamma}}{n-\gamma} \bigg|_0^R \right]\notag\\
=& \frac{\lambda^2 R}{n} + \frac{\mu_1 \lambda R^{1-\beta}}{n-\beta} + \frac{\mu_2 R^{1-\gamma}}{n-\gamma}\notag
\end{align}
under the assumption $n > n-\beta > 0$ and $n-\gamma > 0$. Taking $R = r_c$, we obtain
\begin{align*}
\frac{y_1}{y_0}&=\left| \frac{\phi_\lambda'(1/\lambda)}{\phi_\lambda(1/\lambda)}+\frac{n-1}{2}\lambda \right|\\
&\leq (\frac1n+\frac{n-1}2)\lambda +\frac{\mu_1}{n-\beta}\lambda^{\beta}+\frac{\mu_2}{n-\gamma}\lambda^{\gamma-1}\\
&\lesssim\lambda
\end{align*}
for sufficient small $\lambda$. When $\beta>n$ or $\gamma>n$, both functions $\frac{r^{n-1}}{1+r^{\beta}}$ and $\frac{r^{n-1}}{1+r^{\gamma}}$ are integrable on $[0,\infty)$, and their integrals are bounded by constants. Consequently, in estimate\eqref{ineq:phi}, the last two terms are bounded by $\lambda\mu_1R^{1-n}\sim \lambda^n$ and $\mu_2R^{1-n}\sim \lambda^{n-1}$, respectively. In the boardline case $\beta = n$ or $\gamma=n$, the integral yields
$$\lambda R^{1-n}\ln(1+R^n)=\frac{\ln(1+R^n)}{R^n}\sim \frac1{1+R^n}\sim\lambda^n$$
In both situation, the first term $\la^2R/n\sim \la$ still dominates, and the upper bound for the quotient remains of order $\lambda$.

We now introduce the representation of $y$ for $r>r_c$ in the form
\[ y = y_0\exp\left( \lambda \int_{r_c}^r d\tau + \int_{r_c}^r \mu(\tau) \, d\tau\right)=\exp(\int_{r_c}^r \mu(\tau)d\tau-1)y_0e^{\lambda r}, \]
where $\mu$ is the unknown function. A direct computation yields
\begin{align*}
 y' = y(\lambda + \mu(r)),\quad
y'' = y'(\lambda + \mu(r)) + y \cdot \mu'(r)= y\left[\lambda^2 + \mu^2 + 2\lambda \mu + \mu'\right].
\end{align*}
Substituting these expressions into \eqref{eq:y}, we find that $\mu$ satisfies the Riccati-type equation
\begin{equation}\label{eq:mu}\mu^2+2\lambda\mu+\mu'+G(r)=0,\end{equation}
which immediately implies
\begin{equation}\label{ineq:mu}
\mu'+2\lambda\mu\leq-G(r)
\end{equation}
since $\mu^2\geq0$. We proceed to estimate $\mu(r) = {y'}/{y} - \lambda$ for $r>r_c$. At the boundary $r=r_c$, we have
\[
\mu(r_c) = \frac{y_1}{y_0} - \lambda < \frac{y_1}{y_0}\lesssim \lambda.
\]
Multiplying \eqref{ineq:mu} by the integral factor $ \exp(\int_{r_c}^r 2\lambda \, d\tau)$, we obtain
\[
\left[ \exp\left(\int_{r_c}^r 2\lambda \, d\tau\right) \mu(r) \right]' \leq -G(r) \cdot \exp\left(\int_{r_c}^r 2\lambda \, d\tau\right).
\]
Integrating over interval $[r_c,r]$ gives
\begin{align}
\mu(r) &\leq e^{2\lambda(r_c-r)} \left( \mu(r_c) - \int_{r_c}^r G(s) e^{2\lambda(s-r_c)} \, ds \right)\notag\\
&\leq \mu(r_c) + \int_{r_c}^r|G(s)| \, ds \lesssim \lambda +\int_{r_c}^r|G(s)| \, ds \label{mu(r)}
\end{align}
since \( r_c < s < r \) and $\lambda>0$. We then estimate the integration of $G$.
\begin{eqnarray*}
\int_{r_c}^r|G(s)| \, ds &\leqslant& \int_{r_c}^{r}\left| \frac{(n-1)(n-3)}{4s^2}\right|ds + \lambda \int_{r_c}^{r} \frac{\mu_1}{1+s^\beta} \, ds +\int_{r_c}^{r} \frac{\mu_2}{1+s^\gamma} \, ds
\\
&\leqslant&- \frac{|(n-1)(n-3)|}{4s} \bigg|_{r_c}^{\infty}  + \lambda \mu_1\frac{ s^{1-\beta}}{1-\beta} \bigg|_{r_c}^{\infty}+\mu_2 \frac{s^{1-\gamma}}{1-\gamma} \bigg|_{r_c}^{\infty} \\
&=& \frac{|(n-1)(n-3)|}{4} \lambda + \frac{\lambda \mu_1 \lambda^{\beta-1}}{\beta-1} + \frac{\mu_2 \lambda^{\gamma-1}}{\gamma-1}
\end{eqnarray*}
as $\beta > 1, \gamma > 2$. More precisely, for sufficiently small $\lambda<1$ and $n\neq1,3$, we have
\[\int_{r_c}^r|G(s)| \, ds\leqslant \frac{|(n-1)(n-3)|}{4} \lambda + \frac{\mu_1 \lambda^{\beta}}{\beta-1} + \frac{\mu_2 \lambda^{\gamma-1}}{\gamma-1}\lesssim \lambda .\]
In the case $n=3$, the first term vanishes, while the remaining two terms are bounded by \(\lambda^{\min\{\beta,\gamma-1\}}\), which is less than $\la$. Hence, combining the upperbound of $\mu(r_c)$, we always have $ \mu(r)\lesssim \lambda$.

On the other hand, the lower bound of $\mu(r)$ is given directly from definition $\mu(r) = {y'}/{y} - \lambda>-\la$. Hence, we conclude that
\[
-\lambda < \mu(r) \lesssim\lambda \quad \text{for}\quad r > r_c.
\]
Furthermore, from \eqref{eq:mu} we have
\[
2\lambda\mu = -G - \mu^2 - \mu' \leq -G - \mu',
\]
so integrating over $[r_c,r]$ yields,
\begin{align}
 \int_{r_c}^r \mu(s) \, ds \leq \dfrac{ -\int_{r_c}^r G \, ds - \mu(r)+\mu(r_c)}{2\la}\lesssim1\label{mu_up}
\end{align}
where the last bound uses the $L^1$-estimate of $G$.
To obtain a lower bound for the integral of $\mu$, we rewrite \eqref{eq:mu} as
$-\mu = \frac{\mu'+G}{\mu+2\lambda}$. Then, by integration by parts, we get
\begin{align}
- \int_{r_c}^r \mu(\tau) \, d\tau &= \int_{r_c}^r \frac{\mu'}{\mu + 2\lambda} \, d\tau + \int_{r_c}^r \frac{G}{\mu + 2\lambda} \, d\tau\notag\\
&= \frac{\mu}{\mu+2\lambda} \Bigg|_{r_c}^{r} + \int_{r_c}^{r} \frac{\mu \mu' +(\mu + 2\lambda) G}{(\mu+2\lambda)^2} \, d\tau\notag\\
&= \frac{\mu}{\mu+2\lambda} \Bigg|_{r_c}^{r} + \int_{r_c}^{r} \frac{\mu(\mu'+G)}{(\mu+2\lambda)^2} \, d\tau + \int_{r_c}^{r} \frac{2\lambda G}{(\mu+2\lambda)^2} \, d\tau\notag\\
&\leq \frac{\mu}{\mu+2\lambda} \Bigg|_{r_c}^{r} + \int_{r_c}^{r} \frac{2\lambda G}{(\mu+2\lambda)^2} \, d\tau\notag\\
&\leq \frac{\mu(r)}{\mu(r)+2\lambda}+\frac2\lambda\|G\|_{L^1}\lesssim1\label{mu_low}.
\end{align}
The last two inequalities are due to positivity of $\mu+2\lambda$ and $\mu+2\lambda>\lambda>0$.

Combing \eqref{mu_up} and \eqref{mu_low}, we obtain uniform control on the integral of the integral of $\mu$ over $[r_c,r]$. Consequently, we conclude that
$$
\phi_\lambda(r)=r^{-\frac{n-1}2}y=r^{-\frac{n-1}2}y_0e^{\lambda r}e^{\int_{r_c}^r\mu(s)ds-1}
\backsimeq(\lambda r)^{-\frac{n-1}2}e^{\lambda r}\ ,$$
which completes the proof.
\end{proof}

\begin{rem}
Lemma \ref{exterior} is essential for the critical exponent $p=p_S(n)$. We emphasize the above framework remains valid for the critical damping case \(D(r)=\dfrac{\mu_1}{1+r}\). In fact, we may proceed in exactly the same manner up to the Riccati-type equation \eqref{eq:mu}. We then set
\[
\tilde{\mu}=\mu-\frac{\mu_1}{2(1+r)},
\]
and a direct computation yields
\[
\tilde{\mu}^2+2(\lambda+\frac{\mu_1}{4(1+r)})\tilde{\mu}+\tilde{\mu}'+\tilde{G}=0,
\]
where
\[
\tilde{G}
=\frac{\mu_1^2}{4(1+r)^2}
-\frac{\mu_1}{2(1+r)^2}
-\frac{(n-1)(n-3)}{4r^2}
-\frac{\mu_2}{(1+r)^\gamma}.
\]
Proceeding with \(\tilde{\mu}\) as above, we may find that
\[
\left|\int_{r_c}^r \tilde{\mu}(\tau)\,d\tau\right|\lesssim 1,
\]
and hence establish
\[
\phi_\lambda(r)\backsimeq (\lambda r)^{-\frac{n-\mu_1-1}{2}}e^{\lambda r}.
\]
The exponent shift comes from the multiplier $$\exp(\int_{r_c}^r\frac{\mu_1}{2(1+s)}ds)=(\frac{\la+\la r}{1+\la})^{\frac{\mu_1}2}\backsimeq(\la r)^{\frac{\mu_1}2}$$ for small $\la$.
\end{rem}
\subsection{Asymptotic of $\phi_\lambda(r)$ inside $B_c$}
In order to show the blow-up result for critical exponent, we also need the uniform positive low bound estimate of $\phi_{\lambda}$ inside the critical ball. We first introduce the lemma concerning the solution of limiting equation of $\eqref{eq:test_para}$ as $\lambda$ goes $0$.

\begin{lem}\label{phi_0}
Assume $V=V(r)=\mathcal{O}(<r>^{-\gamma})$ for some $\gamma>2$ is a nonnegative radial potential and is nontrivial in ball $B_{R_0}$.
Then the limiting equation
\begin{equation}\label{limiting_h0}
\Delta \phi_0-V{\phi_0}=0,~ x\in \R^n
\end{equation}
subject to the boundary condition $\phi_0|_{\partial B_{R_0}}=1$ has a unique positive radial solution $\phi_0(r)$ which is increasing outside the ball $B_{R_0}$. Moreover,
for $n\geq3$, $\phi_0$ is bounded and admits the asymptotic expansion
\begin{equation}\phi_0(r) = c_1+c_2r^{2-n}+o(1) \quad \text{as}\quad r\ \rightarrow \infty, \label{asy:h3}\end{equation}
where $c_1>0$.
For $n=2$, the solution behaves like
\begin{equation}\phi_0(r)=\tilde{c}_1+\tilde{c}_2\ln(r)+o(1)\quad \text{as}\quad r\ \rightarrow \infty . \label{asy:h2}\end{equation}
with $\tilde c_2>0$.
\end{lem}
\begin{proof}
The lemma is essentially same as Lemma 2.2 in \cite{LaiLiuTuWang2023}. However, we provide here a direct proof that also establishes the precise asymptotic behavior of the solution. The existence and uniqueness of solution inside the ball $B_{R_0}$ are derived from Lax-Milgram theorem and the maximum principle. Also, the maximum principle ensures $0<\phi_0<1$ in $B_{R_0}$. Moreover, the solution is radial due to the rotation invariance of equation and the boundary condition.

Rewriting \eqref{limiting_h0} in the radial variable, we get the divergence form for $\phi_0(r)$
$$(r^{n-1}\phi_0'(r))'= r^{n-1}V(r)\phi_0(r).$$
Integrating over $[0,r]$ yields
$$\phi_0'(r) = r^{1-n}\int_0^rs^{n-1}V(s)\phi_0(s)ds,$$
which is strictly positive for any $r> R_0$. We may take $(\phi_0(R_0),\phi_0'(R_0))$ as the initial condition for the ODE and thus obtain the global unique solution $\phi_0(r)$. Obviously from the integral, $\phi_0(r)$ is increasing outside the ball $B_{R_0}$ and nondecreasing inside the ball $B_{R_0}$.

Let $r=e^s$ and set $h(s)=\phi_0(e^s)$. Then we derive the elliptic equation for $h(s)$
\begin{equation}\label{eq-0} h''+(n-2)h'-e^{2s}V(e^s)h=0.\end{equation}
Let $Y=(h,\ h_s)^T$. The above equation can be written in following form:
$$\frac d{ds}Y = (\Lambda+R(s))Y,\quad \Lambda =\left(\begin{array}{cc}0 & 1 \\ 0 &  2-n\end{array}\right) \quad R=\left(\begin{array}{cc}0 & 0 \\ e^{2s}V(e^s) & 0\end{array}\right).$$
Since $e^{2s}V(e^{s})=O(e^{(2-\gamma)s})$ and $\gamma>2$, we have $R\in L^1[0,\infty)$. Hence, by the classical Levinson theorem, for $n\geq3$, the two linearly independent solutions possess the asymptotic forms
$$Y_{+} = \left(\begin{array}{c}1+o(1)\\o(1)\end{array}\right),\ Y_{-} = \left(\begin{array}{c}1+o(1)\\2-n+o(1)\end{array}\right)e^{(2-n)s},$$
which imply
$$h(s) = c_1+c_2 e^{(2-n)s}+o(1).$$
Returning to the variable $r$ gives the expansion \eqref{asy:h3}. Moreover, $c_1 = \lim_{r\rightarrow\infty}\phi_0(r)>0$,
and consequently $0<\phi_0(0)\leq \phi_0(r)\leq c_1$, so $\phi_0$ is bounded.

It remains to consider the case $n=2$. Here we could use the extended version of Levinson's theorem (thanks to Hartman-Wintner), as the perturbation $R$ is exponentially decaying.
Alternatively, we give a direct proof.
Let $q(s)=e^{2s}V(e^s)=\mathcal{O}(e^{(2-\gamma)s})$ with $\gamma>2$,
\[
  h_1(s)=1+\int_s^\infty (t-s)q(t)h_1(t)\,dt\ ,~
  h_2(s)=s+\int_s^\infty (t-s)q(t)h_2(t)\,dt\ .
\]
By contraction mapping,
it is clear that $h_j$ exists for $r\gg 1$. In addition,
 $h_j$ are independent solutions to the equation \eqref{eq-0}.
Concerning the asymptotic behavior, we get
$$h_1=1+\mathcal{O}(e^{(2-\gamma)s}),\ h_1'=\mathcal{O}(e^{(2-\gamma)s}),\
h_2=s+\mathcal{O}(e^{(2-\gamma)s}),\ h_2'=1+\mathcal{O}(se^{(2-\gamma)s}) \  .
$$

Thus $h(s) = \tilde{c}_1 h_1+\tilde{c}_2 h_2$ for some $(\tilde{c}_1, \tilde{c}_2)$. Recall that $\phi_0'(r)>0$ for $r\gg 1$, we conclude $\tilde{c}_2>0$, which gives the logarithmic-type asymptotic \eqref{asy:h2}.
\end{proof}

\begin{rem}\label{c_0}
In view of the asymptotic expansions \eqref{asy:h3}, \eqref{asy:h2}, together with the uniform lower bound $\phi_0(r)>\phi_0(0)>0$, we obtain the following asymptotic equivalences as $r\rightarrow\infty$:
$$\phi_0(r)\thicksim 1 \quad (n\geq3)\quad \text{and} \quad \phi_0(r)\thicksim \ln(r+2)\quad (n=2).$$
 Moreover, there exits a constant $r_0>2$, depending only on the solution $\phi_0$, such that for all $r>r_0$, the following strict two-sided bounds hold,
\begin{equation}\frac1{1.1}{c_1} <\phi_0(r)<1.1 c_1 \quad (n\geq3)\quad \text{and} \quad \frac1{1.1}{\tilde{c}_2} <\frac{\phi_0(r)}{\ln(r)}<1.1 \tilde{c}_2 \quad (n=2).\label{two_side}\end{equation}
\end{rem}
We then turn to the asymptotic for solution $\phi_\lambda(r)$ inside the ball $B_c$ by utilizing the asymptotic of $\phi_0$.
\begin{lem}\label{equivlance}
Let $\gamma>2$. Then, inside the critical ball $B_c$, the solution of \eqref{eq:test_radial} is comparable to $\phi_0(r)/\phi_0(1/\lambda)$. More precisely,
\begin{equation}
\phi_\lambda(r)\simeq \frac{\phi_0(r)}{\phi_0(1/\lambda)}
\label{inside}
\end{equation}
for $\lambda<\lambda_1$, where $\lambda_1>0$ is a sufficiently small constant.
\end{lem}

\begin{proof}
Set $u_\lambda(r)=\phi_0(r)/\phi_0(1/\lambda)$ and let $w=\phi_\lambda-u_\lambda$. Then
\[
Lw=w''+\frac{n-1}{r}w'-V(r)w=(\lambda^2+\lambda D)\phi_\lambda>0.
\]
Since $w(1/\lambda)=0$, the maximum principle yields $w\le0$ inside $B_c$. Hence
\[
\phi_\lambda(r)<u_\lambda(r)\qquad\text{for } r<1/\lambda.
\]

For the reverse inequality, set
\[
\psi_\lambda(r)=\frac{\phi_\lambda(r)}{u_\lambda(r)}
=\phi_0(1/\lambda)\frac{\phi_\lambda(r)}{\phi_0(r)}.
\]
We aim to show that $\psi_\lambda$ admits a uniform positive lower bound for small $\lambda$. Clearly, $\psi_\lambda(1/\lambda)=1$, and a direct computation gives
\begin{equation}\label{psi}
(r^{n-1}u^2\psi')'=r^{n-1}(\lambda^2+\lambda D(s))u^2\psi,
\end{equation}
where we suppress the subscript $\lambda$ for notational simplicity. Integrating \eqref{psi} over $[0,r]$, we have
\[
\psi'(r)=\frac{1}{r^{n-1}u^2}\int_0^r s^{n-1}(\lambda^2+\lambda D(s))\psi(s)u^2(s)\,ds.
\]
Integrating again over $[r,1/\lambda]$, we obtain
\begin{align*}
\frac{\phi_\lambda(r)}{u_\lambda(r)}=\psi_\lambda(r)
&=1-\int_r^{1/\lambda}\frac{1}{t^{n-1}u_\lambda^2(t)}\,dt
\int_0^t(\lambda^2+\lambda D)\psi_\lambda(s)s^{n-1}u_\lambda^2(s)\,ds\\
&\ge 1-I(\Omega),
\end{align*}
where
\begin{equation}\label{I(Omega)}
I(\Omega):=\iint_{\Omega}(\lambda^2+\lambda D(s))\psi_\lambda(s)
\frac{s^{n-1}\phi_0^2(s)}{t^{n-1}\phi_0^2(t)}\,ds\,dt,
\end{equation}
with $\Omega:=\{0<s<t,\ 0<t<1/\lambda\}$.
We claim that $I(\Omega)$ admits a uniform upper bound $c_0<1$ for small $\lambda$, so that
\[
\phi_\lambda(r)>(1-c_0)u_\lambda(r).
\]
Combining this with the upper bound, we deduce the desired equivalence in $B_c$ for small $\lambda<\lambda_1$.

It remains to prove
\begin{equation}\label{upper_bound}
\sup_{\lambda<\lambda_1}I(\Omega)<c_0<1.
\end{equation}
We split the domain $\Omega$ into three regions:
\[
\Omega_1=\{0<s<r_0,\ s<t<r_0\},\quad
\Omega_2=\{0<s<r_0,\ r_0<t<1/\lambda\},\quad
\Omega_3=\{r_0<s<t,\ r_0<t<1/\lambda\},
\]
where $r_0$ is the constant from Remark \ref{c_0} depending only on $\phi_0$. Then
\[
I(\Omega)=I(\Omega_1)+I(\Omega_2)+I(\Omega_3).
\]
Noting that $\psi_\lambda(s)<1$ for $s<1/\lambda$ and $s^{n-1}\phi_0^2(s)<t^{n-1}\phi_0^2(t)$ for $s<t$, we have
\[
I(\Omega_1)\le \iint_{\Omega_1}(\lambda^2+\lambda D)\,ds\,dt=O(\lambda).
\]
To bound $I(\Omega_2)$ and $I(\Omega_3)$, we use the two-sided bounds \eqref{two_side}. We first consider $n\ge3$.
\begin{align*}
I(\Omega_2)&\le \frac{1.1^2}{c_1^2}\int_{r_0}^{1/\lambda}\frac{1}{t^{n-1}}\,dt
\int_0^{r_0}(\lambda^2+\lambda D)s^{n-1}\phi_0^2(s)\,ds\\
&=\frac{1.1^2}{c_1^2}\left(\frac{r_0^{2-n}-\lambda^{n-2}}{n-2}\right)O(\lambda)=O(\lambda)
\end{align*}
and
\begin{align*}
I(\Omega_3)&\le 1.1^4\int_{r_0}^{1/\lambda}\,dt\int_{r_0}^t(\lambda^2+\lambda D)\left(\frac{s}{t}\right)^{n-1}\,ds\\
&=1.1^4\int_{r_0}^{1/\lambda}(\lambda^2+\lambda D(s))s^{n-1}\,ds\int_s^{1/\lambda}\frac{1}{t^{n-1}}\,dt\\
&=\frac{1.1^4}{n-2}\int_{r_0}^{1/\lambda}(\lambda^2+\lambda D)s^{n-1}(s^{2-n}-\lambda^{n-2})\,ds\\
&\le \frac{1.1^4}{n-2}\left(\frac{1}{2}-\frac{1}{2}\lambda^2r_0^2
+\frac{r_0^{2-\beta}}{\beta-2}\lambda-\frac{1}{\beta-2}\lambda^{\beta-1}\right)
=\frac{1.1^4}{2(n-2)}+O(\lambda).
\end{align*}
Combining these estimates, we have shown that
\[
I(\Omega)<\frac{1.1^4}{2(n-2)}+\delta=:c_0<1
\]
for sufficiently small $\lambda<\lambda_1$.

Next, we consider $n=2$. The estimates for $I(\Omega_2)$ and $I(\Omega_3)$ are as follows.
\begin{align*}
I(\Omega_2)&\le \frac{1.1^2}{\tilde c_2^2}\int_{r_0}^{1/\lambda}\frac{1}{\ln^2t\,t}\,dt
\int_0^{r_0}(\lambda^2+\lambda D)s\phi_0^2(s)\,ds\\
&=\frac{1.1^2}{\tilde c_2^2}\left(\frac{1}{\ln r_0}+\frac{1}{\ln\lambda}\right)O(\lambda)=O(\lambda)
\end{align*}
and
\begin{align*}
I(\Omega_3)&\le 1.1^4\int_{r_0}^{1/\lambda}(\lambda^2+\lambda D)\ln^2s\,s\,ds
\int_s^{1/\lambda}\ln^{-2}t\,d\ln t\\
&=1.1^4\int_{r_0}^{1/\lambda}(\lambda^2+\lambda D)\ln^2s\,s
\left(\frac{1}{\ln s}-\frac{1}{\ln(1/\lambda)}\right)ds\\
&=1.1^4\int_{r_0\lambda}^1\left(1+\frac{D(\theta/\lambda)}{\lambda}\right)\theta
\left(\ln(\theta/\lambda)+\frac{\ln^2(\theta/\lambda)}{\ln\lambda}\right)d\theta,
\end{align*}
where we have made the change of variables $\theta=\lambda s$ in the last identity. Since
\[
\frac{D(\theta/\lambda)}{\lambda}
=\frac{\mu_1\lambda^{\beta-1}}{\lambda^\beta+\theta^\beta}
<\frac{\mu_1}{\theta^\beta}\lambda^{\beta-1},\qquad \beta>1,
\]
we have
\begin{align*}
I(\Omega_3)&\le 1.1^4\int_{r_0\lambda}^{1}
\left(1+\frac{\mu_1}{\theta^\beta}\lambda^{\beta-1}\right)\theta
\left(\frac{\ln^2\theta}{\ln\lambda}-\ln\theta\right)d\theta\\
&\le 1.1^4\int_0^1\theta
\left(\frac{\ln^2\theta}{\ln\lambda}-\ln\theta\right)d\theta
+1.1^4\int_{r_0\lambda}^1\mu_1\left(\frac{\lambda}{\theta}\right)^{\beta-1}
\left(\frac{\ln^2\theta}{\ln\lambda}-\ln\theta\right)d\theta\\
&=\frac{1.1^4}{4}+O\left(\frac{1}{\ln\lambda}\right)
+1.1^4\mu_1\left(\frac{\lambda^{\beta-1}}{\ln\lambda}\int_{r_0\lambda}^1\theta^{1-\beta}\ln^2\theta\,d\theta
-\lambda^{\beta-1}\int_{r_0\lambda}^1\theta^{1-\beta}\ln\theta\,d\theta\right).
\end{align*}
It remains to treat the integral terms in the last bracket, which arise from the damping term. If $\beta=2$, these terms reduce to
\[
\frac{\lambda}{\ln\lambda}\int_{r_0\lambda}^1\frac{\ln^2\theta}{\theta}\,d\theta
-\lambda\int_{r_0\lambda}^1\frac{\ln\theta}{\theta}\,d\theta
=-\frac{\lambda\ln^3(r_0\lambda)}{3\ln\lambda}
-\frac{\lambda\ln^2(r_0\lambda)}{2}
=O(\lambda\ln^2\lambda).
\]
If $\beta\ne2$, a direct computation gives
\begin{align*}
\int_{\eta}^{1}\theta^{1-\beta}\ln\theta\,d\theta
&=\frac{\eta^{2-\beta}}{\beta-2}\ln\eta
+\frac{\eta^{2-\beta}}{(2-\beta)^2}
-\frac{1}{(2-\beta)^2}
\end{align*}
and
\begin{align*}
\int_{\eta}^{1}\theta^{1-\beta}\ln^2\theta\,d\theta
&=\frac{\eta^{2-\beta}}{\beta-2}\ln^2\eta
+\frac{2}{(\beta-2)^2}\eta^{2-\beta}\ln\eta
-\frac{2}{(2-\beta)^3}\eta^{2-\beta}
+\frac{2}{(2-\beta)^3}.
\end{align*}
Plugging in $\eta=r_0\lambda$, we find that these terms are at least $O(\lambda\ln\lambda)$ for $\beta>2$, or $O(\lambda^{\beta-1})$ for $1<\beta<2$, both of which tend to $0$ as $\lambda\to0$. Summarizing the above computations, we conclude that
\[
I(\Omega)\le \frac{1.1^4}{4}+\delta=:c_0<1
\]
for sufficiently small $\lambda$. Hence \eqref{upper_bound} is established for all $n\ge2$.
\end{proof}

\begin{rem}Combing the results of Lemma \ref{phi_0} and Lemma \ref{equivlance}, we obtain the following equivalence inside the critical ball $B_c$:
 $$\phi_\lambda(r)\backsimeq  \frac{\phi_0(r)}{\phi_0(\frac1\lambda)}\backsimeq\left\{\begin{array}{lcl}\frac{\phi_0(r)}{\ln(\frac1\lambda)}\ &\text{for}& n=2\\
 \phi_0(r)\ &\text{for}& n\geq3\end{array}\right.$$
 for small $\lambda\leq\lambda_0:=\min(\lambda_1,\lambda_2,1)$. Here $\lambda_2:=\frac1{r_0}$, so that $\frac1\lambda>r_0$. Furthermore, the monotonicity of $\phi_0$ implies $\phi_\la(r)$ admits a uniform upper bound in the critical ball.
\end{rem}

\section{Non-radial solution $\phi_\lambda(x)$ and test function $b_q(x,t)$}
In this section, we derive the asymptotic behavior of solution $\phi_\la(x)$ without imposing the radial symmetry assumption on the coefficients $D$ and $V$. The main idea is to employ the comparison principle together with a monotone exhaustion argument. Based on the resulting estimates, we then construct the test function $b_q(x,t)$ and establish its asymptotic properties.

\subsection{Asymptotic of $\phi_\lambda(x)$}
For \(n\geq2\) and \(0<\la\le\la_0\), we define
\begin{equation}
\label{eq:Phi}
\Phi_\la(x)
:=
\frac1{|\mathbb S^{n-1}|}
\int_{\mathbb S^{n-1}}e^{\la x\cdot\omega}\,d\omega,
\end{equation}
which satisfies \(
\Delta\Phi_\la=\la^2\Phi_\la
\) and $\Phi_\la>0$ in $\R^n$. Moreover, for \(0<\la\le\la_0\), the following uniform estimate holds:
\begin{equation}
\label{eq:Phi-est}
\Phi_\la(x)
\asymp
(1+\la|x|)^{-\frac{n-1}{2}}e^{\la|x|}
\end{equation}
uniformly for \(x\in\R^n\). We are now in a position to state the following theorem concerning the existence and estimates of the solution to \eqref{eq:test_para}.
\begin{thm}
\label{thm:elliptic-global}
Let \(n\ge2\) and assume \eqref{eq:DV-main-ass}. Then for every \(0<\lambda\le\lambda_0\) there exists a positive weak solution
\[
\phi_\lambda\in H^1_{\rm loc}(\mathbb R^n)\cap W^{2,s}_{\rm loc}(\mathbb R^n)
\quad\text{for every }1<s<\infty
\]
satisfying \eqref{eq:test_para}. Moreover, there exists a sufficiently small constant \(c_0>0\), independent of \(\lambda\) and \(x\), such that
\begin{equation}
\label{eq:phi-barrier}
c_0\tilde{\phi}_\lambda(|x|)\le\phi_\lambda(x)\le\Phi_\lambda(x),
\qquad x\in\mathbb R^n,
\end{equation}
where \(\tilde{\phi}_\lambda\) is the corresponding radial barrier solution introduced in Section 2. The family \(\{\phi_\lambda\}\) may be chosen so that \(\lambda\mapsto\phi_\lambda(x)\) is measurable for every fixed \(x\).
\end{thm}

\begin{proof}
Set
\[
Q_\lambda(x):=\lambda^2+\lambda D(x)+V(x),
\qquad
Q_\lambda^*(x):=\lambda^2+\lambda D_*(|x|)+V_*(|x|),
\]
and \(L_\lambda:=-\Delta+Q_\lambda(x)\), \(L_\lambda^*:=-\Delta+Q_\lambda^*(x)\). Since
\[
L_\lambda\Phi_\lambda=(\lambda D+V)\Phi_\lambda\ge0,
\]
\(\Phi_\lambda\) is a positive supersolution. On the other hand, since \(L_\lambda^*\tilde{\phi}_\lambda=0\) and \(Q_\lambda\le Q_\lambda^*\), we have
\[
L_\lambda \tilde{\phi}_\lambda
=(L_\lambda-L_\lambda^*)\tilde{\phi}_\lambda
=(Q_\lambda-Q_\lambda^*)\tilde{\phi}_\lambda\le0,
\]
so \(\tilde{\phi}_\lambda\) is a positive subsolution.

Recalling the uniform upper bound of \(\tilde{\phi}_\lambda(r)\) inside the critical ball and its asymptotic behavior outside, we can choose a sufficiently small constant \(c_0>0\), independent of \(\lambda\), such that
\[
c_0\tilde{\phi}_\lambda(|x|)<\frac12\Phi_\lambda(|x|)
\]
for all \(x\in\mathbb R^n\). For each integer \(N\ge1\), let \(\phi_{\lambda,N}\) be the weak solution of the truncated Dirichlet problem
\begin{equation}
\label{eq:exhaust}
\begin{cases}
L_\lambda\phi_{\lambda,N}=0,&x\in B_N,\\
\phi_{\lambda,N}=\Phi_\lambda,&x\in\partial B_N.
\end{cases}
\end{equation}
Writing \(\phi_{\lambda,N}=\Phi_\lambda+w\), we see that \(w\) satisfies
\[
-\Delta w+Q_\lambda w=-(\lambda D+V)\Phi_\lambda.
\]
Since \(Q_\lambda\ge\lambda^2>0\), \(D,V\in L^\infty\), and \(\Phi_\lambda\) is smooth and bounded on \(B_N\), the Lax--Milgram theorem yields a unique \(w\in H_0^1(B_N)\), and hence a unique weak solution \(\phi_{\lambda,N}\in H^1(B_N)\). Applying the maximum principle on \(B_N\) gives
\begin{equation}
\label{eq:barrier-N}
c_0\tilde{\phi}_\lambda(|x|)\le\phi_{\lambda,N}(x)\le\Phi_\lambda(x)
\quad\text{a.e. in }B_N.
\end{equation}
If \(N_2>N_1\), then on \(\partial B_{N_1}\subset B_{N_2}\),
\[
\phi_{\lambda,N_2}\le\Phi_\lambda=\phi_{\lambda,N_1},
\]
so the comparison principle yields monotonicity in \(N\):
\[
\phi_{\lambda,N_2}\le\phi_{\lambda,N_1}
\qquad\text{in }B_{N_1}.
\]
Thus the pointwise limit
\[
\phi_\lambda(x):=\lim_{N\to\infty}\phi_{\lambda,N}(x)
\]
exists and is positive by \eqref{eq:barrier-N}.

We show this pointwise limit is weak solution of \eqref{eq:test_para}.
On every compact \(K\Subset\R^n\), \eqref{eq:barrier-N} gives a
uniform \(L^\infty(K)\) bound. Since \(Q_\la\in L^\infty\), the standard
interior elliptic \(W^{2,s}\)-estimates implies that the limit belongs
to \(W^{2,s}_{\rm loc}\) for every \(1<s<\infty\) and solves
\eqref{eq:test_para} in the sense of distribution. The comparison \eqref{eq:phi-barrier} is then obtained by
passing to the limit.

For fixed \(N\), the solution of \eqref{eq:exhaust} depends
continuously on \(\la\) in \(H^1(B_N)\), and interior elliptic
estimates upgrade this to local uniform continuity. Since
\(\ph_\la\) is the monotone limit along the integer exhaustion, the
family can be chosen measurable in \(\la\).
\end{proof}

\begin{rem} Bounded measurability of $Q_\la$ is sufficient for the weak Dirichlet problems,
the weak maximum principle, and the local \(W^{2,s}\)-regularity used above. Hence, we don't require any regularity for $D$ and $V$.
\end{rem}

\begin{rem} Combining the results of Lemma \ref{exterior}, Lemma \ref{equivlance} and asymptotic \eqref{eq:Phi-est}, the two side bound \eqref{eq:phi-barrier} implies $\phi_\la(x)$ shares same asymptotic as $\tilde{\phi}_\la(|x|)$, i.e, for $n\geq2$, outside the critcal ball $\la |x|>1$,
$$\phi_\la(x)\backsimeq (\la|x|)^{-\frac{n-1}{2}}e^{\la|x|};$$
and inside the critical ball $\la |x|<1$,
$$\phi_\la(x)\sim 1 \quad (n\geq3) \quad \text{and} \quad \frac{\tilde{\phi}_0(|x|)}{\tilde{\phi}_0(1/\la)}\lesssim \phi_\la(x)\lesssim 1 \quad (n=2). $$
\end{rem}

\subsection{Asymptotic of $b_q(x,t)$}
Once the asymptotic behavior of \(\phi_\lambda(x)\) has been established, we can construct the test function for the critical case in the standard manner. For \(q>0\), let
\begin{equation}\label{def-bq}
b_q(t,x)
:=
\int_0^{\lambda_0} e^{-\lambda t}\phi_\lambda(x)\,\lambda^{q-1}\,d\lambda.
\end{equation}
Clearly, \(b_q\) satisfies the adjoint equation
\begin{equation}
\left(\partial_t^2-\Delta-D(x)\partial_t+V(x)\right)b_q=0
\label{bq-PDE}
\end{equation}
and
\begin{equation}
\partial_t b_q=-b_{q+1}.
\label{bq-der}
\end{equation}
We then derive the asymptotic behavior of \(b_q(x,t)\) from that of \(\phi_\lambda\), which is crucial for our final blow-up result. In this subsection, we establish the asymptotic expansion for \(b_q(x,t)\). By utilizing the behavior of \(\phi_\lambda(x)\) both inside and outside the critical ball \(B_c\), we obtain the following lemma.
\begin{lem}
For \(n\ge3\), the test function \(b_q(x,t)\) satisfies the following estimates for \(|x|<t+R\):
\begin{align}
b_q(x,t) &\gtrsim (t+R)^{-q}\phi_0(|x|),\qquad \text{if}\quad q>0,\label{bq-lower-3}\\
b_q(x,t) &\lesssim (t+R)^{-q},\qquad \text{if}\quad 0<q<\frac{n-1}{2},\label{bq-upper1}\\
b_q(x,t) &\lesssim (t+R)^{-\frac{n-1}{2}}\,(t+R+1-|x|)^{\frac{n-1}{2}-q},
\qquad \text{if}\quad q>\frac{n-1}{2}.\label{bq-upper2}
\end{align}
For \(n=2\), the upper bounds are the same as \eqref{bq-upper1} and \eqref{bq-upper2}, while the lower bound reads
\begin{equation}
b_q(x,t) \gtrsim \frac{\phi_0(|x|)}{\ln(t+R)}(t+R)^{-q},\qquad \text{if}\quad q>0.\label{bq-lower-2}
\end{equation}
An additional logarithmic factor appears due to the logarithmic-type lower bound of \(\phi_\lambda\) inside the critical ball.
\end{lem}

\begin{proof}
The proof of this lemma is essentially the same as that in \cite{LaiLiuTuWang2023}. For the reader's convenience, we provide the details here. Set \(r=|x|\). To be precise, we write
\begin{align*}
b_q(x, t) &= \int_0^{\lambda_0} \phi_\lambda(x) e^{-\lambda t} \lambda^{q-1} \, d\lambda\\
&= \int_0^{\min(1/r,\lambda_0)} \phi_\lambda(x) e^{-\lambda t} \lambda^{q-1} \, d\lambda
 + \int_{\min(1/r,\lambda_0)}^{\lambda_0}\phi_\lambda(x) e^{-\lambda t} \lambda^{q-1} \, d\lambda\\
&= A + B.
\end{align*}
Here \(A\) and \(B\) denote the integrals inside and outside the critical ball \(B_c\), respectively. The term \(B\) vanishes when \(r<1/\lambda_0\).

We first consider \(n\ge3\). Since \(\lambda<\lambda_0\), we may replace \(\phi_\lambda(x)\) by \(\phi_0(|x|)/\phi_0(1/\lambda)\) in \(A\), so that
\[
A \simeq \phi_0(r)\int_0^{\min(1/r,\lambda_0)}\frac{e^{-\lambda t} \lambda^{q-1}}{\phi_0(1/\lambda)} \, d\lambda.
\]
As \(r<1/\lambda\), we have \(\phi_0(r)<\phi_0(1/\lambda)\), hence for any \(q>0\),
\begin{equation}\label{A_upper}
A\lesssim \int_0^{\min(1/r,\lambda_0)}e^{-\lambda t} \lambda^{q-1} \, d\lambda
\lesssim (t+R)^{-q}\int_0^{\infty}e^{-\theta}\theta^{q-1} \, d\theta
\simeq (t+R)^{-q}.
\end{equation}
For the lower bound of \(A\), since \(\phi_0(1/\lambda)\simeq c_1\) for \(\lambda<\lambda_0\), we have for any \(q>0\),
\begin{align}
A&\gtrsim \phi_0(r)\int_0^{\min(1/r,\lambda_0)}e^{-\lambda t} \lambda^{q-1} \, d\lambda\notag \\
&= \phi_0(r)(t+R)^{-q}\int_0^{\min((t+R)/r,\lambda_0(t+R))}e^{-\theta} \theta^{q-1} \, d\theta\notag\\
&> \phi_0(r)(t+R)^{-q}\int_0^{\min(1,\lambda_0R)}e^{-\theta}\theta^{q-1} \, d\theta
\simeq \phi_0(r)(t+R)^{-q}.\label{A_lower_3}
\end{align}
We then turn to the upper bound of \(B\). Since \(r\ge1/\lambda>1/\lambda_0\), we apply \eqref{outside} for \(\phi_\lambda(r)\).

\(\text{(1)}\quad r < \frac{t+R}{2}\) for any \(q>0\).
\begin{align}
B &= \int_{1/r}^{\lambda_0} \phi_\lambda(x) e^{-\lambda t} \lambda^{q-1} \, d\lambda
\approx \int_{1/r}^{\lambda_0} e^{\lambda r} (\lambda r)^{-\frac{n-1}{2}} e^{-\lambda t} \lambda^{q-1} \, d\lambda\notag\\
&\lesssim \int_{1/r}^{\lambda_0} e^{\lambda (r-t)} \lambda^{q-1} \, d\lambda
< \int_{1/r}^{\lambda_0} e^{\lambda ((t+R)/2-t)} \lambda^{q-1} \, d\lambda\notag\\
&\lesssim e^{\lambda_0R} \int_{1/r}^{\lambda_0} e^{-\lambda (t+R)/2} \lambda^{q-1} \, d\lambda
\simeq \left(\frac{t+R}{2}\right)^{-q} \int_{(t+R)/(2r)}^{(t+R)\lambda_0/2} e^{-\theta} \theta^{q-1} \, d\theta\notag\\
&\lesssim (t+R)^{-q} \int_0^{\infty} e^{-\theta} \theta^{q-1} \, d\theta\notag\\
&\le (t+R)^{-q}.\label{B_upper_1}
\end{align}

\(\text{(2)}\quad \frac{t+R}{2} \le r \le t+R\) for \(0<q<\frac{n-1}{2}\).
\begin{align}
B &\le \int_{1/r}^{\lambda_0} e^{-\lambda t} (\lambda r)^{-\frac{n-1}{2}} e^{\lambda r} \lambda^{q-1} \, d\lambda
\le \int_{1/r}^{\lambda_0} e^{\lambda R} (\lambda r)^{-\frac{n-1}{2}} \lambda^{q-1} \, d\lambda\notag\\
&\lesssim r^{-q}\int_1^{r\lambda_0} \theta^{-\frac{n-1}{2}} \theta^{q-1} \, d\theta
= r^{-q} \int_1^{r\lambda_0} \theta^{q-\frac{n-1}{2}-1} \, d\theta\notag\\
&\lesssim(t+R)^{-q},\label{B_upper_2}
\end{align}
since \(r\simeq t+R\) and \(0<q<\frac{n-1}{2}\).

\(\text{(3)}\quad \frac{t+R}{2} \le r \le t+R\) for \(q>\frac{n-1}{2}\).
\begin{align}
B &\le \int_{1/r}^{\lambda_0} e^{-\lambda t} (\lambda r)^{-\frac{n-1}{2}} e^{\lambda r} \lambda^{q-1} \, d\lambda
= r^{-\frac{n-1}{2}}\int_{1/r}^{\lambda_0} e^{-\lambda t + \lambda r} \lambda^{q-\frac{n-1}{2}-1} \, d\lambda\notag\\
&\lesssim (t+R)^{-\frac{n-1}{2}} \int_{1/r}^{\infty} e^{-\lambda (t+R+1-r)} \lambda^{q-\frac{n-1}{2}-1} \, d\lambda\notag\\
&\le (t+R)^{-\frac{n-1}{2}}(t+R+1-r)^{\frac{n-1}{2}-q}.\label{B_upper_3}
\end{align}

We conclude the result for \(n\ge3\) from the above analysis. The lower bound \eqref{bq-lower-3} follows from \eqref{A_lower_3} and the fact that \(b_q(x,t)\ge A\). The upper bound \eqref{bq-upper1} for \(q<\frac{n-1}{2}\) is established from \eqref{A_upper}, \eqref{B_upper_1}, and \eqref{B_upper_2}. The upper bound \eqref{bq-upper2} follows from \eqref{B_upper_3} for \(q>\frac{n-1}{2}\). It is easy to see that the upper bound provided by \(B\) is dominant.

Finally, we consider the case \(n=2\). It is easy to check that all the estimates hold as in the case \(n\ge3\), except for the lower bound of \(A\), which is modified due to the different asymptotic behavior of \(\phi_0\) inside the critical ball \(B_c\). Indeed, we have \(\phi_0(1/\lambda)\simeq \ln(1/\lambda)\) for \(\lambda<\lambda_0\), so that
\begin{align}
A &\simeq\phi_0(r)\int_0^{\min(1/r,\lambda_0)} \frac{e^{-\lambda (t+R)}}{\ln(1/\lambda)} \, \lambda^{q-1}d\lambda \notag \\
&\ge (t+R)^{-q} \phi_0(r)\int_0^{\min((t+R)/r,\lambda_0(t+R))}\frac{ e^{-\theta}\theta^{q-1}}{\ln(t+R)-\ln\theta} \, d\theta \notag \\
&\ge (t+R)^{-q} \phi_0(r)\int_{1/(t+R)}^{\min(1,\lambda_0R)} \frac{e^{-\theta}\theta^{q-1}}{\ln(t+R)-\ln\theta} \, d\theta\notag\\
&\gtrsim \frac{\phi_0(r)}{\ln(t+R)}(t+R)^{-q}\label{A_lower_2}
\end{align}
for \(t>1/(\lambda_0R)-R\). If \(t\le 1/(\lambda_0R)-R\),
\begin{align*}
A&\gtrsim (t+R)^{-q} \phi_0(r)\int_0^{\min(1,\lambda_0R)} \frac{e^{-\theta}\theta^{q-1}}{-2\ln\theta} \, d\theta\\
&\gtrsim (t+R)^{-q}\phi_0(r)I_0\frac{\ln R}{\ln(t+R)}
\gtrsim \frac{\phi_0(r)}{\ln(t+R)}(t+R)^{-q}.
\end{align*}
This provides the lower bound \eqref{bq-lower-2} for \(n=2\).
\end{proof}

\section{The proof of main theorem: Critical case with potential exponent $\gamma>2$}
In this section, we present the proof of the main theorem. Owing to the careful preparation of the test function, the potential term vanishes directly after taking integration by parts. The argument follows the standard approach stated in \cite{LaiLiuTuWang2023}. Indeed, it suffices to verify inequality $(6.4)$ and $(6.15)$ therein, and further combining which with the asymptotic estimates \eqref{bq-lower-3} through \eqref{bq-lower-2} completes the proof. For the reader's convenience, we give the proof. For case $n=2$, some additional treatment is required due to the different lower bound \eqref{bq-lower-2} and we therefore provide more details.

Let $\eta$ be the smooth cut-off function
\[
\eta(t) = \left\{
\begin{array}{ll}
1, \, \quad &\text{for}\quad t \leq \frac{1}{2}\\
\text{decreasing},\, \quad &\text{for}\quad \frac{1}{2} < t < 1\\
0, \, \quad &\text{for}\quad t \geq 1
\end{array}
\right.\]
and set \( \eta_T(t) = \eta(\frac{t}T) \). Denoting $\Phi(t,x)=e^{-\lambda t}\phi_{\lambda}(x)$ and $\Psi(t,x)=\eta^{2p'}_T(t)\Phi(t,x)$,
we start the proof from testing \eqref{eq:PDE0} by $\Psi(t,x)$. Taking integration by parts to derive
\begin{align*}
&\int_0^T\int_{\mathbb{R}^n}|u|^p\Psi dxdt-\int_{\mathbb{R}^n}\left(u_t(t,x)+D(x)u(t,x)\right)\Psi dx\bigg|_{0}^T+\int_{\mathbb{R}^n}u(t,x)\partial_t\Psi dx\bigg|_{0}^T\\&=\int_0^T\int_{B_{t+R}}u(\partial_t^2-\Delta-D\partial_t+V)\Psi dxdt.
\end{align*}
Noting that $\eta_T(0)=1$, $\eta_T(T)=0$, $\eta'_T(0)=\eta'_T(T)=0$ and $\partial_t\Phi=-\Phi$, we deduce that
\begin{align}
&\int_0^T \int_{\mathbb{R}^n} |u|^p \Psi dx dt + \varepsilon \int_{\mathbb{R}^n} (g(x) + (\lambda + D(x)) f(x)) \Phi(0,x) dx\notag\\
&= \int_0^T \int_{B_{t+R}} u(\partial_t^2 - \Delta - D\partial_t + V) \Psi dx dt\notag\\
&= \int_0^T \int_{B_{t+R}} u\bigl(\partial_t^2 (\eta_T^{2p'}) + 2\partial_t (\eta_T^{2p'}) \partial_t - D\partial_t (\eta_T^{2p'}) \bigr) \Phi
dx dt\notag\\
&= \int_0^T \int_{B_{t+R}} u\bigl(\partial_t^2 (\eta_T^{2p'}) - 2\lambda \partial_t (\eta_T^{2p'}) - D\partial_t (\eta_T^{2p'}) \bigr) \Phi dx dt\notag\\
&\leq C \int_{\frac{T}2}^T \int_{B_{t+R}} |u| \, \eta_T^{2(p'-1)} \bigl(T^{-2} + (2\lambda + |D|) T^{-1}\bigr) \Phi dx dt=\text{RHS}.\label{mid1}
\end{align}
We estimate the integral in the righthandside. Applying H\"{o}lder inequality gives
\begin{align}
\text{RHS}\leq \frac{C}T \left(\int_{\frac{T}2}^T \int_{B_{t+R}} |u|^p \, \eta_T^{2p'}dxdt\right)^{\frac1p}\notag
\left(\int_{\frac{T}2}^T \int_{B_{t+R}} \Phi^{p'} dx dt\right)^{\frac1{p'}}\notag\\
\leq CT^{-1+[1+(n-1)(1-\frac{p'}2)]/p'} \left(\int_{\frac{T}2}^T \int_{B_{t+R}} |u|^p \, \eta_T^{2p'}dxdt\right)^{\frac1p} \label{mid2}
\end{align}
where we estimate the integral concerning $\Phi$ by utilizing Lemma 2.1 and Lemma 4.1 in \cite{LaiLiuTuWang2023}, as
\begin{align*}
\int_{\frac{T}2}^T \int_{B_{t+R}} \Phi^{p'} dx dt&\leq \int_{\frac T2}^{T}\int_0^{t+R}e^{-p'\lambda t}(\lambda r)^{-\frac{n-1}2p'}e^{p'\lambda r}r^{n-1}drdt\\
&\leq \int_{\frac T2}^T(t+R)^{(n-1)(1-\frac{p'}2)}dt\approx T^{1+(n-1)(1-\frac{p'}2)}.\end{align*}
Hence, plugging \eqref{mid2} into \eqref{mid1}, we get
\begin{equation}\label{6.4}
\int_{T/2}^T \int_{\mathbb{R}^n} \eta_T^{2p'} |u|^p dx dt\geq C_2^p(f,g)  \varepsilon^p T^{p}T^{-\frac{p}{p'}[1+(n-1)(1-\frac{p'}2)]} =C\varepsilon^p T^{n-\frac{n-1}{2}p}
\end{equation}
where $$C_2(f,g)=\int_{\mathbb{R}^n} (g(x) + (\lambda + D(x)) f(x)) \phi_\lambda(x) dx.$$

Second, we test \eqref{eq:PDE0} by $\Psi(x,t)= b_q(x,t)\eta_T^{2p'}(t)$ where $q=\frac{n-1}2-\frac1p$ and $p=p_S(n)$
\begin{align*}
&\int_0^T \int_{\mathbb{R}^n} b_q\eta_T^{2p'}|u|^{p} dx dt
= \int_0^T \int_{\mathbb{R}^n} \left( \partial_t^2 - \Delta + D(x) \partial_t +V(x)\right)u b_q \eta_T^{2p'} dx dt \\
&\lesssim \int_0^T \int_{\mathbb{R}^n} |u| \left( |2\partial_t b_q \partial_t \eta_T^{2p'}| + |b_q \partial_t^2 \eta_T^{2p'}| + |D(x) b_q \partial_t \eta_T^{2p'}| \right) dx dt \\
&\leq \left( \int_{\frac{T}{2}}^T \int_{\mathbb{R}^n} |u|^p \Psi dx dt \right)^{\frac{1}{p}} \left( \int_{\frac{T}{2}}^T \int_{B_{t+R}} b_q^{-\frac{p'}{p}} \left[ (T^{-1} + D)b_q + b_{q+1} \right]^{p'} T^{-p'} dx dt \right)^{\frac{1}{p'}}.
\end{align*}
Note the integral in second bracket can be bounded by:
\[\left(\int_{\frac{T}{2}}^T \int_{B_{t+R}} (T^{-1} + D)^{p'} T^{-p'} b_q dx dt\right)^{\frac1{p'}}+\left(\int_{\frac{T}{2}}^T \int_{B_{t+R}} b_{q+1}^{p'} b_q^{-\frac{p'}{p}} T^{-p'} dx dt\right)^{\frac1{p'}}=I^{\frac1{p'}}+II^{\frac1{p'}}.\]
As \( D \leq \mu_1(1+r)^{-\beta} \) and \eqref{bq-upper1}, we have for \( p = p_S(n) \),
\begin{align*}
I\leq \int_{\frac{T}{2}}^T \int_{0}^{t+R} (1+r)^{n-1-p'} T^{-p'} (t+R)^{-q}dr dt\leq T^{-2p'-q+n+1}=1
\end{align*}
As $q+1=\frac{n+1}{2}-\frac1p>\frac{n-1}2$, we have
\(b_{q+1} \lesssim(t+R)^{-\frac{n-1}{2}} (t+R+1-r)^{\frac{n-1}{2}-(q+1)}\) and
\[\quad b_q(x,t)\gtrsim\left\{\begin{array}{ll} \frac{\phi_0(r)}{\ln(t+R)}(t+R)^{-q}\quad&\text{for}\quad n=2,\\
(t+R)^{-q}&\text{for}\quad n\geq3.
\end{array}\right.\]
Plugging those into $II$. For $n\geq3$, we have
\begin{align*}b_{q+1}^{p'} b_q^{-\frac{p'}{p}} T^{-p'} &\leq (t+R)^{-\frac{n-1}{2}p'} (t+R+1-r)^{(\frac{n-1}{2}-(q+1))p'}(t+R)^{q\frac{p'}{p}}T^{-p'}\\
&=(t+R)^{-\frac{n-1}2p'+(\frac{n-1}{2}-\frac1p)(p'-1)}T^{-p'}(t+R+1-r)^{-1}\\
&=T^{-p'-\frac{n-1}2+\frac{1-p'}p}(t+R+1-r)^{-1}.
\end{align*}
Hence, for $p=p_S(n)$,
\begin{align*}II&\leq \int_{\frac{T}{2}}^T \int_0^{t+R} T^{\frac{n-1}2-p'+\frac{1-p'}p}(t+R+1-r)^{-1} dr dt \\
&= \int_{\frac{T}{2}}^T \int_0^{t+R} T^{-1}(t+R+1-r)^{-1} dr dt \leq \ln(T). \end{align*}
For $n=2$, we have additional factor $(\frac{\ln(t+R)}{\phi_0(r)})^{\frac{p'}{p}}$ from $b_q^{-\frac{p'}{p}}$,
\begin{align*}
II&\leq \int^{T}_{\frac T2}\int_0^{t+R}T^{-1}(t+R+1-t)^{-1}(\frac{\ln(t+R)}{\phi_0(r)})^{\frac{p'}{p}}drdt
\end{align*}
Since $\phi_0(r)\backsim\ln(r+2)$, we can see $\frac{\ln(t+R)}{\phi_0(r)}\backsim1$ if $r>\frac{t+R}2$. The integral estimation of this part is reduced to same result as $n\geq3$. We only focus on the part $r<\frac{t+R}2$,
\begin{align*}&\int^{T}_{\frac T2}\int_0^{\frac{t+R}{2}}T^{-1}(t+R+1-r)^{-1}(\frac{\ln(t+R)}{\ln(r+2)})^{\frac{p'}{p}}drdt \\
&\leq \int^{T}_{\frac T2}\int_0^{\frac{t+R}{2}}T^{-1}(r+2)^{-1}(\frac{\ln(t+R)}{\ln(r+2)})^{\frac{p'}{p}}drdt\\
&= T^{-1}\int^{T}_{\frac T2}\ln^{\frac{p'}{p}}(t+R)dt\int_0^{\frac{t+R}{2}}(r+2)^{-1}\ln^{-\frac{p'}p}(r+2)dr\leq \ln(T).
\end{align*}
Hence, for both $n=2$ and $n\geq3$, we have
\begin{equation}\label{X(T)}
X(T) \leq (\ln(T))^{\frac1{p'}} Z(T)^{\frac1{p}}\end{equation}
where
\[X(T) \triangleq \int_0^T \int_{\mathbb{R}^n} |u|^p b_q \eta_T^{2p'} dx dt\quad \text{and} \quad
Z(T) \triangleq \int_{T/2}^T \int_{\mathbb{R}^n} |u|^p b_q \eta_T^{2p'} dx dt.
\]

Now, we just follow the same argument in \cite{LaiLiuTuWang2023} to obtain the blow-up result for both $n\geq3$ and $n=2$ with same lifespan estimate. Let $Z=TY'$ and $Y(2)=0$ then
$$Y(M)=\int_2^MZ(T)T^{-1}dT.$$
Direct calculation gives
\[Y(M)\leq \ln(2)X(M),\]
plugging which into \eqref{X(T)}gives
$$Y(T)\lesssim(\ln T)^{1/{p'}}(TY'(T))^{1/p}\quad \text{or} \quad Y'(T)Y^{-p}\gtrsim T^{-1}(\ln T)^{1-p},\quad \forall T\in [2,T_{\epsilon}).$$
Integrating the above inequality from $T_1$ to $T_2$ gives for $p=p_S(n)>1$
\begin{equation}\label{Y_1}
Y(T_1)^{1-p}\geq Y(T_2)^{1-p}+c\int_{\ln T_1}^{\ln T_2}\tau^{1-p}d\tau \quad \forall 2<T_1<T_2<T_\varepsilon.\end{equation}

Recalling \eqref{6.4} and plugging the lower bound of $b_q$, i.e, \eqref{bq-lower-3} and \eqref{bq-lower-2} for $n\geq3$ and $n=2$ respectively, we get the lower bound of $Z(t)$
\begin{align*}
TY'=Z(T)= \int_{\frac T2}^T\int_{\mathbb{R}^n}|u|^p \eta^{2p'}_Tb_qdxdt\geq  \int_{\frac T2}^T\int_{|x|>\frac{t+R}2}|u|^p \eta^{2p'}_Tb_qdxdt\gtrsim \varepsilon^p T^{n-\frac{n-1}2p-q}=\varepsilon^p,
\end{align*}
noticing that the lower bound of $b_q$ is same for $n=2$ and $n\geq3$ near the light cone $|x|>\frac{t+R}2$.
Integrating $Y'\gtrsim \frac{\varepsilon^p}T$ from $2$ to $T_1$ gives
\begin{equation}\label{Y_2}
Y(T_1)>Y(2)+c \varepsilon^p(\ln T_1-\ln2)\gtrsim \varepsilon^p\ln T_1,\quad \forall T_1\in (4, T_\varepsilon).
\end{equation}

Combining \eqref{Y_1} and \eqref{Y_2}, we have
$$c\int_{\ln T_1}^{\ln T_2}\tau^{1-p}d\tau \leq Y(T_1)^{1-p} \leq \varepsilon^{p(1-p)}(\ln T_1)^{1-p}.$$
Let $T_1=\sqrt{T_2}$ and $T_2\rightarrow T_\varepsilon$, we find
$$\ln T_\varepsilon \leq C_p\varepsilon^{p(1-p)}$$
which gives the desired lifespan stated in main theorem.

\paragraph*{\bf Acknowledgment}

The first author was supported by JSPS KAKENHI Grant Numbers 25K24910 and 23K03174.
The third author was supported by NSFC 12671293.

\appendix

\section{Local existence and finite propagation speed}
\label{sec:appendix-local}

We consider
\begin{equation}
\label{eq:appendix-PDE}
\begin{cases}
u_{tt}-\Delta u+D(x)u_t+V(x)u=|u|^p,
& (t,x)\in (0,T)\times\mathbb R^n,\\
u(0,x)=\varepsilon f(x),\qquad
u_t(0,x)=\varepsilon g(x),
& x\in\mathbb R^n.
\end{cases}
\end{equation}
Assume that \(D\) and \(V\) are real-valued measurable functions satisfying
\begin{equation}
\label{eq:appendix-DV}
|D(x)|
\le
\frac{\mu_1}{(1+|x|)^\beta},
\qquad
|V(x)|
\le
\frac{\mu_2}{(1+|x|)^\gamma},
\end{equation}
where
\(
\mu_1,\  \mu_2\ge0,
\
\beta>1,
\
\gamma>2.
\)

We emphasize no sign condition on \(D\) or \(V\) is needed in this appendix.

\begin{prop}
\label{prop:appendix-local}
Let
\(
\varepsilon>0,~
f\in H^1(\mathbb R^n),
~
g\in L^2(\mathbb R^n).
\)
Assume
\[ 1<p<\infty\quad \text{if}\quad n=2\quad \text{and}\quad 1<p\le\frac{n}{n-2}\quad \text{if}\quad n\ge3.\]
Then there exists a unique energy solution of \eqref{eq:appendix-PDE}
\[
u\in
C([0,T];H^1(\mathbb R^n))
\cap
C^1([0,T];L^2(\mathbb R^n))
\]
with
\[
T
=
T\left(
\varepsilon,
\|f\|_{H^1},
\|g\|_{L^2},
\|D\|_{L^\infty},
\|V\|_{L^\infty}
\right)
>0.
\]
If \(T_\varepsilon\in(0,\infty]\) denotes the maximal existence time,
then either
\(
T_\varepsilon=\infty
\)
or
\[
\limsup_{t\uparrow T_\varepsilon}
\left(
\|u(t)\|_{H^1}
+
\|u_t(t)\|_{L^2}
\right)
=
\infty.
\]
If, in addition,
\[
\operatorname{supp}f\cup\operatorname{supp}g
\subset \overline{B}_{R_0}
\]
for some \(R_0>0\), then the solution has finite propagation speed:
\[
\operatorname{supp}u(t,\cdot)
\subset \overline{B}_{R_0+t},
\qquad 0\le t<T_\varepsilon.
\]
\end{prop}

\begin{proof}
Let
\[
X_T
:=
C([0,T];H^1(\mathbb R^n))
\cap
C^1([0,T];L^2(\mathbb R^n))
\]
with norm
\[
\|v\|_{X_T}
:=
\sup_{0\le t\le T}
\left(
\|v(t)\|_{H^1}
+
\|v_t(t)\|_{L^2}
\right).
\]
For \(M>0\), set
\[
X_T(M)
:=
\left\{
v\in X_T:
\|v\|_{X_T}\le M
\right\}.
\]

For \(v\in X_T(M)\), let \(\mathcal Mv=u\) be the solution of
\begin{equation}
\label{eq:appendix-map}
\begin{cases}
u_{tt}-\Delta u
=
|v|^p-D(x)v_t-V(x)v,\\
u(0)=\varepsilon f,\qquad
u_t(0)=\varepsilon g.
\end{cases}
\end{equation}
The energy estimate for the inhomogeneous wave equation gives
\begin{align}
\|\mathcal Mv(t)\|_{H^1}
+
\|\partial_t\mathcal Mv(t)\|_{L^2}
&\le
C\varepsilon
\left(
\|f\|_{H^1}
+
\|g\|_{L^2}
\right)
\notag\\
&\quad
+
C\int_0^t
\left\|
|v(s)|^p
-
Dv_t(s)
-
Vv(s)
\right\|_{L^2}
\,ds.
\label{eq:appendix-energy-map}
\end{align}

\noindent
We estimate the right-hand side. Since \(D,V\in L^\infty\),
\[
\|Dv_t(s)\|_{L^2}
\le
\|D\|_{L^\infty}\|v_t(s)\|_{L^2}
\]
and
\[
\|Vv(s)\|_{L^2}
\le
\|V\|_{L^\infty}\|v(s)\|_{L^2}.
\]
Thus
\begin{equation}
\label{eq:appendix-linear-lower}
\|Dv_t(s)+Vv(s)\|_{L^2}
\le
C_{D,V}M,
\end{equation}
where
\[
C_{D,V}
:=
\|D\|_{L^\infty}
+
\|V\|_{L^\infty}.
\]

Moreover,
\[
\||v(s)|^p\|_{L^2}
=
\|v(s)\|_{L^{2p}}^p.
\]
By Sobolev embedding,
\[
H^1(\mathbb R^n)
\hookrightarrow
L^{2p}(\mathbb R^n),
\]
and hence
\begin{equation}
\label{eq:appendix-nonlinearity}
\||v(s)|^p\|_{L^2}
\le
C\|v(s)\|_{H^1}^p
\le
CM^p.
\end{equation}

Combining
\eqref{eq:appendix-energy-map},
\eqref{eq:appendix-linear-lower},
and
\eqref{eq:appendix-nonlinearity},
we obtain
\[
\|\mathcal Mv\|_{X_T}
\le
C\varepsilon
\left(
\|f\|_{H^1}
+
\|g\|_{L^2}
\right)
+
CT
\left(
C_{D,V}M+M^p
\right).
\]

Choose
\[
M
=
2C\varepsilon
\left(
\|f\|_{H^1}
+
\|g\|_{L^2}
\right)
+
1.
\]
Taking \(T>0\) sufficiently small, we obtain
\[
\mathcal M:X_T(M)\longrightarrow X_T(M).
\]

Let
\[
v_1,v_2\in X_T(M),
\qquad
u_j=\mathcal Mv_j,
\]
and set
\[
w=u_1-u_2,
\qquad
z=v_1-v_2.
\]
Then
\[
w_{tt}-\Delta w
=
|v_1|^p-|v_2|^p-Dz_t-Vz,
\qquad
w(0)=w_t(0)=0.
\]
Therefore
\begin{align}
\|w\|_{X_T}
&\le
C\int_0^T
\||v_1|^p-|v_2|^p\|_{L^2}
\,ds
\notag\\
&\quad
+
C\int_0^T
\left(
\|Dz_t\|_{L^2}
+
\|Vz\|_{L^2}
\right)
\,ds.
\end{align}

Since
\[
\bigl||a|^p-|b|^p\bigr|
\le
C_p
\left(
|a|^{p-1}+|b|^{p-1}
\right)|a-b|,
\]
H\"older's inequality and Sobolev embedding yield
\[
\begin{aligned}
\||v_1|^p-|v_2|^p\|_{L^2}
&\le
C
\left(
\|v_1\|_{L^{2p}}^{p-1}
+
\|v_2\|_{L^{2p}}^{p-1}
\right)
\|v_1-v_2\|_{L^{2p}}\\
&\le
CM^{p-1}\|z\|_{X_T}.
\end{aligned}
\]
Also,
\[
\|Dz_t\|_{L^2}
+
\|Vz\|_{L^2}
\le
C_{D,V}\|z\|_{X_T}.
\]
Hence
\[
\|\mathcal Mv_1-\mathcal Mv_2\|_{X_T}
\le
CT
\left(
C_{D,V}+M^{p-1}
\right)
\|v_1-v_2\|_{X_T}.
\]
Taking \(T>0\) smaller if necessary, \(\mathcal M\) is a contraction.
Banach's fixed-point theorem gives a unique local energy solution.

It remains to prove the finite propagation property when the initial data
are compactly supported. We use the Picard sequence associated with the
same contraction mapping. Set
\[
u_0=0,
\qquad
u_{m+1}=\mathcal M u_m,
\qquad m\ge0.
\]
We claim that
\begin{equation}
\label{eq:appendix-Picard-support}
\operatorname{supp}u_m(t,\cdot)
\subset \overline{B}_{R_0+t},
\qquad 0\le t\le T,
\end{equation}
for every \(m\ge0\). This is immediate for \(u_0=0\). Assume that
\eqref{eq:appendix-Picard-support} holds for \(u_m\). By the definition
of \(\mathcal M\), \(u_{m+1}\) solves
\[
\begin{cases}
\partial_t^2u_{m+1}-\Delta u_{m+1}
=|u_m|^p-D(x)\partial_tu_m-V(x)u_m,\\
u_{m+1}(0)=\varepsilon f,\qquad
\partial_tu_{m+1}(0)=\varepsilon g.
\end{cases}
\]
Since differentiation with respect to time does not enlarge spatial
support, the induction hypothesis implies
\[
\operatorname{supp}
\bigl(|u_m(t)|^p-D\partial_tu_m(t)-Vu_m(t)\bigr)
\subset \overline{B}_{R_0+t}.
\]
The finite propagation speed for the free inhomogeneous wave equation
therefore yields
\[
\operatorname{supp}u_{m+1}(t,\cdot)
\subset \overline{B}_{R_0+t},
\qquad 0\le t\le T.
\]
Thus \eqref{eq:appendix-Picard-support} follows by induction.

Since the Picard sequence converges to the fixed point in \(X_T\),
\[
u_m\longrightarrow u
\quad\text{in }X_T,
\]
and hence, for every fixed \(t\in[0,T]\),
\(u_m(t)\to u(t)\) in \(L^2(\mathbb R^n)\). Because every
\(u_m(t)\) vanishes almost everywhere on
\(\mathbb R^n\setminus\overline{B}_{R_0+t}\), the limit \(u(t)\)
also vanishes there. Consequently,
\[
\operatorname{supp}u(t,\cdot)
\subset \overline{B}_{R_0+t},
\qquad 0\le t\le T.
\]

The continuation criterion follows by restarting the same local
argument at any time \(t_0<T_\varepsilon\). The same restart argument,
together with the already established support property up to \(t_0\),
extends the finite propagation property to every
\(0\le t<T_\varepsilon\).
\end{proof}

\begin{rem}
\label{rem:appendix-Strauss-Sobolev}
For \(n\ge3\), the above argument requires
\[
p\le\frac{n}{n-2}.
\]
Moreover,
\[
p_S(3)<3,
\qquad
p_S(4)=2,
\qquad
p_S(n)>\frac{n}{n-2}
\quad (n\ge5).
\]
Thus the argument covers \(p=p_S(n)\) for \(n=2,3,4\). However, for \(n\ge5\), a Strichartz-based local theory is required.
\end{rem}

\bibliographystyle{plain}

\end{CJK*}

\end{document}